\documentclass{article}
\usepackage{amsthm}
\usepackage{amssymb}
\usepackage{amsmath}
\usepackage{amsfonts}
\usepackage{color}
\usepackage{indentfirst} 
\usepackage{graphicx}
\usepackage{color}
\usepackage[square, comma, sort&compress, numbers]{natbib}

\newtheorem{lemma}{Lemma}[section]
\newtheorem{theorem}{Theorem}[section]

\newtheorem{corollary}{Corollary}[section]

\def\b1{\mbox{\boldmath $1$}}

\newenvironment{demo*}{\vspace{3mm}\noindent{\bf Proof.}}{\hfill $\Box$ \vspace{3mm}}

\begin{document}
\title{\bf \Large {Conditional tail risk expectations for location-scale mixture of elliptical distributions}}
{\color{red}{\author{
\normalsize{Baishuai  Zuo,\; Chuancun Yin*}\\
{\normalsize\it  School of Statistics, Qufu Normal University, Qufu, Shandong 273165, P. R. China}\\
*Corresponding author. Email address:  ccyin@qfnu.edu.cn (Chuancun Yin)}}}

\maketitle
\vskip0.01cm
\noindent{\large {\bf Abstract}}  We present general results on the univariate tail conditional expectation (TCE) and multivariate tail conditional expectation for location-scale mixture of elliptical distributions. Examples include the location-scale mixture of normal distributions, location-scale mixture of Student-$t$ distributions, location-scale mixture of Logistic distributions and location-scale mixture of Laplace distributions.  We also consider portfolio risk decomposition with TCE for location-scale mixture of elliptical distributions.



\baselineskip=20pt

\section{Introduction and Motivation}

Tail conditional expectation (TCE), one of important risk measures, is common and practical.
TCE of a random variable $X$ is defined as
$TCE_{X}(x_{q})=E(X|X>x_{q}),$
where $x_{q}$ is a particular value, generally referred to as the $q$-th quantile with
$\overline{F}_{X}(x_{q})=1-q.$
Here $\overline{F}_{X}(x)=1-F_{X}(x)$ is tail distribution function of $X$.
The TCE   has been discussed in many literatures ( see Landsman and Valdez (2003), Ignatieva and Landsman (2015,~2019)).
  Recently, a new type of  multivariate tail conditional expectation (MTCE) was defined by Landsman et al. (2016). It is the following special case when $\boldsymbol{q}=(q,~q,\cdots,q)$.
\begin{align*}
&MTCE_{\boldsymbol{q}}(\mathbf{X})=E\left[\mathbf{X}|\mathbf{X}>VaR_{\boldsymbol{q}}(\mathbf{X})\right]\\
&=E[\mathbf{X}|X_{1}>VaR_{q_{1}}(X_{1}),\cdots,X_{n}>VaR_{q_{n}}(X_{n})],~\boldsymbol{q}=(q_{1},\cdots,q_{n})\in(0,~1)^{n},
\end{align*}
 where $\mathbf{X}=(X_{1},~X_{2},\cdots,X_{n})^{T}$ is an $n\times1$ vector of risks with cumulative distribution function $F_{\mathbf{X}}(\boldsymbol{x})$ and tail function $\overline{F}_{\mathbf{X}}(\boldsymbol{x})$, $$VaR_{\boldsymbol{q}}(\mathbf{X})=(VaR_{q_{1}}(X_{1}),~VaR_{q_{2}}(X_{2}),\cdots,VaR_{q_{n}}(X_{n}))^{T},$$ and $VaR_{q_{k}}(X_{k}),~k=1,~2,\cdots,n$ is the value at risk (VaR) measure of the random variable $X_{k}$, being the $q_{k}$-th quantile of $X_{k}$ (see Landsman et al. (2016)). On the basis of it, Mousavi et al. (2019) study  multivariate tail conditional expectation for scale mixtures of skew-normal distribution.

 Closely related to tail conditional expectation is portfolio risk decomposition with TCE, it's research has experienced a rapid growth in the literature. Portfolio risk decomposition based on TCE for the elliptical distribution was studied in Landsman and Valdez (2003) and extended to the multivariate skew-normal distribution in Vernic (2006). The phase-type distributions and multivariate Gamma distribution were researched in Cai and Li (2005) and Furman and Landsman (2007), respectively. Furthermore, Hashorva and Ratovomirija (2015) considered the capital allocation with TCE for mixed Erlang distributed risks joined by the Sarmanov distribution, and Ignatieva and Landsman (2019) has given the expression of TCE-based allocation for the generalised hyperbolic distribution. Recently, Zuo and Yin (2020) deals with the tail conditional expectation
for univariate generalized skew-elliptical distributions and multivariate tail conditional expectation for generalized skew-elliptical distributions.

The rest of the paper is organized as follows. In Section 2 we define the location-scale mixture of elliptical distributions and establish some properties of it. Furthermore, we give several examples as special cases of it. In Section 3 we derive TCE for univariate cases of mixture of elliptical distributions, and in Section 4, we provide expression of  MTCE for mixture of elliptical distributions. Section 5 offers expression of portfolio risk decomposition with TCE for mixture of elliptical distributions. Section 6 gives concluding remarks.

\section{Mixture of elliptical distributions}
In this section we introduce the location-scale mixture of elliptical (LSME) distributions and some its properties.
Let $\mathbf{Y}\sim LSME_{n}(\boldsymbol{\mu},~\mathbf{\Sigma},~\boldsymbol{\beta},~\Theta,~g_{n})$ be an $n$-dimensional LSME distribution with location parameter $\boldsymbol{\mu}$ and positive definite scale matrix $\mathbf{\Sigma}=(\sigma_{i,j})_{i,j=1}^{n}$, if
\begin{align}\label{(1)}
\mathbf{Y}=\boldsymbol{\mu}+\Theta\boldsymbol{\beta}+\Theta^{\frac{1}{2}}\mathbf{\Sigma}^{\frac{1}{2}}\mathbf{X},
\end{align}
where $\boldsymbol{\beta}\in \mathbb{R}^{n}$, and $\mathbf{X}\sim E_{n}(\boldsymbol{0},~\mathbf{I}_{n},~g_{n}).$ Assume that $\mathbf{X}$ is independent of non-negative scalar random variable $\Theta$. We have
\begin{align}\label{(2)}
\mathbf{Y}|\Theta=\theta\sim E_{n}(\boldsymbol{\mu}+\theta\boldsymbol{\beta},~\theta\mathbf{\Sigma},~g_{n}).
\end{align}
 Here $\mathbf{X}$ is an $n$-dimensional elliptical random vector, and denoted by $\mathbf{X}\sim E_{n}\left(\boldsymbol{\mu},~\boldsymbol{\Sigma},~g_{n}\right)$. If it's probability density function exists, the form will be
\begin{align}\label{(3)}
f_{\boldsymbol{X}}(\boldsymbol{x}):=\frac{1}{\sqrt{|\boldsymbol{\Sigma}|}}g_{n}\left\{\frac{1}{2}(\boldsymbol{x}-\boldsymbol{\mu})^{T}\mathbf{\Sigma}^{-1}(\boldsymbol{x}-\boldsymbol{\mu})\right\},~\boldsymbol{x}\in\mathbb{R}^{n},
\end{align}
 where $\boldsymbol{\mu}$ is an $n\times1$ location vector, $\mathbf{\Sigma}$ is an $n\times n$ scale matrix and $g_{n}(u)$, $u\geq0$, is the density generator of $\mathbf{X}$. This density generator satisfies condition: (see Fang et al. (1990))
 \begin{align}\label{(4)}
 \int_{0}^{\infty}u^{\frac{n}{2}-1}g_{n}(u)\mathrm{d} u<\infty.
 \end{align}
The characteristic function of $\mathbf{X}$ takes the form  $\varphi_{\boldsymbol{X}}(\boldsymbol{t})=\exp\left\{i\boldsymbol{t}^{T}\boldsymbol{\mu}\right\}\psi\left(\frac{1}{2}\boldsymbol{t}^{T}\boldsymbol{\Sigma}\boldsymbol{t}\right),~\boldsymbol{t}\in \mathbb{R}^{n}$, with function $\psi(t):[0,\infty)\rightarrow\mathbb{R},$ called the characteristic generator. Furthermore, the condition
 \begin{align}\label{(5)}
 |\psi'(0)|<\infty,
 \end{align}
 guarantees the existence of the covariance matrix of $\mathbf{X}$ (see Fang et al. (1990)).
 Suppose $\mathbf{A}$ is a $k\times n$ matrix, and $\mathbf{b}$ is a $k\times1$ vector.  Then
\begin{align}\label{(6)}
\mathbf{AX+b}\sim E_{k}\left(\boldsymbol{A\mu+b},~\boldsymbol{A\Sigma A^{T}},~g_{k}\right).
\end{align}

 To express conditional tail risk measures for $n$-dimensional mixture of elliptical distributions
we introduce the cumulative generator $\overline{G}_{n}(u)$ (see Landsman et al.(2018)):
\begin{align}\label{(7)}
\overline{G}_{n}(u)=\int_{u}^{\infty}{g}_{n}(v)\mathrm{d}v.
\end{align}

Let $\mathbf{X}^{\ast}\sim E_{n}(\boldsymbol{\mu},~\boldsymbol{\Sigma},~\overline{G}_{n})$ be an elliptical random vector with generator $\overline{G}_{n}(u)$, whose the density function (if it exists)
\begin{align}\label{(8)}
f_{\boldsymbol{X}^{\ast}}(\boldsymbol{x})=\frac{-1}{\psi'(0)\sqrt{|\boldsymbol{\Sigma}|}}\overline{G}_{n}\left\{\frac{1}{2}(\boldsymbol{x}-\boldsymbol{\mu})^{T}\mathbf{\Sigma}^{-1}(\boldsymbol{x}-\boldsymbol{\mu})\right\},~\boldsymbol{x}\in\mathbb{R}^{n}.
\end{align}

We list some examples of the mixture of elliptical family, including location-scale mixture of normal (LSMN) distributions, location-scale mixture of Student-$t$ (LSMSt) distributions, location-scale mixture of Logistic (LSMLo) distributions and location-scale mixture of Laplace (LSMLa) distributions.\\
$\mathbf{Example~2.1}$ (Mixture of normal distribution). An $n$-dimensional normal random vector $\mathbf{X}$ with location parameter $\boldsymbol{\mu}$ and scale matrix  $\mathbf{\Sigma}$  has density function
\begin{align*}
f_{\boldsymbol{X}}(\boldsymbol{x})=\frac{(2\pi)^{-\frac{n}{2}}}{\sqrt{|\boldsymbol{\Sigma}|}}\exp\left\{-\frac{1}{2}(\boldsymbol{x}-\boldsymbol{\mu})^{T}\mathbf{\Sigma}^{-1}(\boldsymbol{x}-\boldsymbol{\mu})\right\},~\boldsymbol{x}\in\mathbb{R}^{n},
\end{align*}
and denoted by $\mathbf{X}\sim N_{n}\left(\boldsymbol{\mu},~\boldsymbol{\Sigma}\right)$. Therefore, the location-scale mixture of normal random vector \\$\mathbf{Y}\sim LSMN_{n}(\boldsymbol{\mu},~\mathbf{\Sigma},~\boldsymbol{\beta},~\Theta)$ is defined as
\begin{align}\label{(20)}
\mathbf{Y}=\boldsymbol{\mu}+\Theta\boldsymbol{\beta}+\Theta^{\frac{1}{2}}\mathbf{\Sigma}^{\frac{1}{2}}\mathbf{X},
\end{align}
and $\boldsymbol{\mu}$, $\mathbf{\Sigma}$, $\Theta$~and~$\boldsymbol{\beta}$ are the same as in $(\ref{(1)})$. \\
 $\mathbf{Example~2.2}$ (Mixture of student-$t$ distribution). An $n$-dimensional student-$t$ random vector $\mathbf{X}$ with location parameter $\boldsymbol{\mu}$, scale matrix $\mathbf{\Sigma}$ and $m>0$ degrees of freedom  has density function
 \begin{align*}
f_{\boldsymbol{X}}(\boldsymbol{x})=&\frac{c_{n}}{\sqrt{|\mathbf{\Sigma}|}}\left[1+\frac{(\boldsymbol{x}-\boldsymbol{\mu})^{T}\mathbf{\Sigma}^{-1}(\boldsymbol{x}-\boldsymbol{\mu})}{m}\right]^{-\frac{m+n}{2}},~\boldsymbol{x}\in\mathbb{R}^{n},
\end{align*}
 where $c_{n}=\frac{\Gamma\left(\frac{m+n}{2}\right)}{\Gamma(m/2)(m\pi)^{\frac{n}{2}}}$, and denoted by $\mathbf{X}\sim St_{n}\left(\boldsymbol{\mu},~\boldsymbol{\Sigma},~m\right)$. So that the location-scale mixture of student-$t$ random vector $\mathbf{Y}\sim LSMSt_{n}(\boldsymbol{\mu},~\mathbf{\Sigma},~\boldsymbol{\beta},~\Theta)$ satisfies
\begin{align}\label{(21)}
\mathbf{Y}=\boldsymbol{\mu}+\Theta\boldsymbol{\beta}+\Theta^{\frac{1}{2}}\mathbf{\Sigma}^{\frac{1}{2}}\mathbf{X},
\end{align}
and $\boldsymbol{\mu}$, $\mathbf{\Sigma}$, $\Theta$~and~$\boldsymbol{\beta}$ are the same as in $(\ref{(1)})$. \\
$\mathbf{Example~2.3}$ (Mixture of Logistic distribution). Density function of an $n$-dimension Logistic random vector $\mathbf{X}$ with location parameter $\boldsymbol{\mu}$ and scale matrix  $\mathbf{\Sigma}$  can be expressed as
\begin{align*}
f_{\boldsymbol{X}}(\boldsymbol{x})=\frac{c_{n}}{\sqrt{|\mathbf{\Sigma}|}}\frac{\exp\left\{-\frac{1}{2}(\boldsymbol{x}-\boldsymbol{\mu})^{T}\mathbf{\Sigma}^{-1}(\boldsymbol{x}-\boldsymbol{\mu})\right\}}{\left[1+\exp\left\{-\frac{1}{2}(\boldsymbol{x}-\boldsymbol{\mu})^{T}\mathbf{\Sigma}^{-1}(\boldsymbol{x}-\boldsymbol{\mu})\right\}\right]^{2}},~\boldsymbol{x}\in\mathbb{R}^{n},
\end{align*}
where $c_{n}=(2\pi)^{-n/2}\left[\sum_{i=0}^{\infty}(-1)^{i-1}i^{1-n/2}\right]^{-1}$, and denoted by $\mathbf{X}\sim Lo_{n}\left(\boldsymbol{\mu},~\boldsymbol{\Sigma}\right)$.
  The location-scale mixture of Logistic random vector $\mathbf{Y}\sim LSMLo_{n}(\boldsymbol{\mu},~\mathbf{\Sigma},~\boldsymbol{\beta},~\Theta)$ satisfies
\begin{align}\label{(22)}
\mathbf{Y}=\boldsymbol{\mu}+\Theta\boldsymbol{\beta}+\Theta^{\frac{1}{2}}\mathbf{\Sigma}^{\frac{1}{2}}\mathbf{X},
\end{align}
and $\boldsymbol{\mu}$, $\mathbf{\Sigma}$, $\Theta$~and~$\boldsymbol{\beta}$ are the same as in $(\ref{(1)})$. \\
 $\mathbf{Example~2.4.}$ (Mixture of Laplace distribution). Density of  Laplace random vector $\mathbf{X}$ with location parameter $\boldsymbol{\mu}$ and scale matrix  $\mathbf{\Sigma}$  is given by
 \begin{align*}
f_{\boldsymbol{X}}(\boldsymbol{x})=&\frac{c_{n}}{\sqrt{|\mathbf{\Sigma}|}}\exp\left\{-[(\boldsymbol{x}-\boldsymbol{\mu})^{T}\mathbf{\Sigma}^{-1}(\boldsymbol{x}-\boldsymbol{\mu})]^{1/2}\right\},~\boldsymbol{x}\in\mathbb{R}^{n},
\end{align*}
where $c_{n}=\frac{\Gamma(n/2)}{2\pi^{n/2}\Gamma(n)}$, and denoted by $\mathbf{X}\sim La_{n}\left(\boldsymbol{\mu},~\boldsymbol{\Sigma}\right)$. Hence, the location-scale mixture of Laplace random vector $\mathbf{Y}\sim LSMLa_{n}(\boldsymbol{\mu},~\mathbf{\Sigma},~\boldsymbol{\beta},~\Theta)$ is defined as
\begin{align}\label{(23)}
\mathbf{Y}=\boldsymbol{\mu}+\Theta\boldsymbol{\beta}+\Theta^{\frac{1}{2}}\mathbf{\Sigma}^{\frac{1}{2}}\mathbf{X},
\end{align}
and $\boldsymbol{\mu}$, $\mathbf{\Sigma}$, $\Theta$~and~$\boldsymbol{\beta}$ are the same as in $(\ref{(1)})$.
\section{Univariate cases}
\begin{theorem}\label{th.1} Let $Y\sim LSME_{1}(\mu,~\sigma^{2},~\beta,~\Theta,~g_{1})$ be an univariate location-scale mixture of elliptical random variable defined as $(\ref{(1)})$.  We suppose
\begin{align}\label{(9)}
\int_{0}^{\infty}g_{1}(u)\mathrm{d}u<\infty.
\end{align}
 Then
\begin{align}\label{(10)}
 &TCE_{Y}(y_{q})=\mu+E_{\theta}\left[\theta\beta+\delta_{\theta}(\sqrt{\theta}\sigma)^{2}\right],
\end{align}
where
\begin{align*}
 \delta_{\theta}=\frac{\frac{1}{\sqrt{\theta}\sigma}\overline{G}_{1}(\frac{1}{2}z_{q}^{2})}{\overline{F}_{Z}(z_{q})},
\end{align*}
with $Z\sim E_{1}(0,~1,~g_{1})$ and $z_{q}=\frac{y_{q}-\mu-\theta\beta}{\sqrt{\theta}\sigma}$.
\end{theorem}
\noindent Proof. Using definition and tower property of expectations, we obtain
\begin{align*}
TCE_{Y}(y_{q})&=E[Y|Y>y_{q}]\\
&=E_{\Theta}[E(Y|Y>y_{q},\Theta)].
\end{align*}
Since
\begin{align*}
E[Y|Y>y_{q},\Theta=\theta]&=E[(Y|\Theta=\theta)|(Y|\Theta=\theta)>y_{q}]&\\
&=E[M|M>y_{q}]\\
&=TCE_{M}(y_{q}),
\end{align*}
where $M\sim E_{1}(\mu+\theta\beta,~\theta\sigma^{2},~g_{1})$, and the second equality we have used (\ref{(2)}).\\
Using Theorem 1 in Landsman and Valdez (2003), we obtain $(\ref{(10)})$, which completes the proof of Theorem \ref{th.1}.\\
$\mathbf{Remark~3.1.}$ We find that $TCE_{Y|\Theta}(y_{q})$ is a special case of Theorem 1 in Landsman and Valdez (2003).
\begin{corollary}\label{co.1}
Let $Y\sim LSMN_{1}(\mu,~\sigma^{2},~\beta,~\Theta)$ be an univariate location-scale mixture of normal random variable defined as $(\ref{(20)})$. Under conditions in Theorem \ref{th.1}, we obtain the TCE for location-scale mixture of normal distributions. Its' form is the same as $(\ref{(10)})$,
where
\begin{align*}
 \delta_{\theta}=\frac{\frac{1}{\sqrt{\theta}\sigma}\phi_{1}(\frac{1}{2}z_{q}^{2})}{1-\Phi_{1}(z_{q})}.
\end{align*}
 Additionally, $\phi_{1}(\cdot)$ and $\Phi_{1}(\cdot)$ denote the density and distribution functions of normal distributions,  respectively.
\end{corollary}
\noindent Proof. Letting the density generator $\overline{G}_{1}(u)=g_{1}(u)=\phi_{1}(u)=(2\pi)^{-\frac{1}{2}}e^{-u}$ in Theorem \ref{th.1}, we directly obtain our result. This completes the proof of Corollary \ref{co.1}.
\begin{corollary}\label{co.2}
Let $Y\sim LSMSt_{1}(\mu,~\sigma^{2},~\beta,~\Theta)$ be an univariate location-scale mixture of student-$t$ random variable defined as $(\ref{(21)})$. Under conditions in Theorem \ref{th.1}, we obtain the TCE for location-scale mixture of Student-$t$ distributions. Its' form is the same as $(\ref{(10)})$,
where
\begin{align*}
 \delta_{\theta}=\frac{\frac{1}{\sqrt{\theta}\sigma}\overline{G}_{1}(\frac{1}{2}z_{q}^{2})}{\overline{F}_{Z}(z_{q})}=\frac{\frac{1}{\sqrt{\theta}\sigma}c_{1}\frac{m}{m-1}(1+\frac{z_{q}^{2}}{m})^{-(m-1)/2}}{\overline{F}_{Z}(z_{q})}=\frac{\frac{1}{\sqrt{\theta}\sigma}t_{m,1}(z_{q};0,1)}{\overline{T}_{m,1}(z_{q};0,1)}.
\end{align*}
In addition, $t_{m,1}(z_{q};0,1)$ and $T_{m,1}(z_{q};0,1)$ are the density and distribution functions of Student-$t$ distributions, respectively (see Landsman et al. (2016)).
\end{corollary}
\noindent Proof. Letting $g_{1}(u)=c_{1}(1+\frac{2u}{m})^{-(m+1)/2}$, $\overline{G}_{1}(u)=c_{1}\frac{m}{m-1}(1+\frac{2u}{m})^{-(m-1)/2}$ and $c_{1}=\frac{\Gamma((m+1)/2)}{\Gamma(m/2)(m\pi)^{\frac{1}{2}}}$ (see Landsman et al. (2016)) in Theorem \ref{th.1}, we immediately obtain our result. This completes the proof of Corollary \ref{co.2}.
\begin{corollary}\label{co.3}
Let $Y\sim LSMLo_{1}(\mu,~\sigma^{2},~\beta,~\Theta)$ be an univariate location-scale mixture of Logistic random variable defined as $(\ref{(22)})$. Under conditions in Theorem \ref{th.1}, we obtain the TCE for location-scale mixture of Logistic distributions. Its' form is the same as $(\ref{(10)})$,
where
\begin{align*}
 \delta_{\theta}=\frac{\frac{1}{\sqrt{\theta}\sigma}\overline{G}_{1}(\frac{1}{2}z_{q}^{2})}{\overline{F}_{Z}(z_{q})}=\frac{\frac{1}{\sqrt{\theta}\sigma}c_{1}\frac{\exp(-\frac{1}{2}z_{q}^{2})}{1+\exp(-\frac{1}{2}z_{q}^{2})}}{\overline{F}_{Z}(z_{q})}=\left[\frac{1}{2[(\sqrt{2\pi})^{-1}+\phi(z_{q})]}\right]\frac{\frac{1}{\sqrt{\theta}\sigma}\phi(z_{q})}{\overline{F}_{z}(z_{q})}.
\end{align*}
In addition, $\phi(\cdot)$ is the density functions of normal distributions (see Landsman and Valdez (2003)).
\end{corollary}
\noindent Proof. Letting $g_{1}(u)=c_{1}\frac{\exp(-u)}{[1+\exp(-u)]^{2}}$, $\overline{G}_{1}(u)=c_{1}\frac{\exp(-u)}{1+\exp(-u)}$ and $c_{1}=\frac{1}{2}$ (see Landsman and Valdez (2003)) in Theorem \ref{th.1}, we directly obtain our result. This completes the proof of Corollary \ref{co.3}.
\begin{corollary}\label{co.4}
Let $Y\sim LSMLa_{1}(\mu,~\sigma^{2},~\beta,~\Theta)$ be an univariate location-scale mixture of Laplace random variable defined as $(\ref{(22)})$. Under conditions in Theorem \ref{th.1}, we obtain the TCE for location-scale mixture of Laplace distributions. Its' form is the same as $(\ref{(10)})$,
where
\begin{align*}
 \delta_{\theta}=\frac{\frac{1}{\sqrt{\theta}\sigma}\overline{G}_{1}(\frac{1}{2}z_{q}^{2})}{\overline{F}_{Z}(z_{q})}=\frac{\frac{1}{\sqrt{\theta}\sigma}c_{1}(1+\sqrt{z_{q}^{2}})\exp(-\sqrt{z_{q}^{2}})}{\overline{F}_{Z}(z_{q})}=\sqrt{2}\left(1+\sqrt{z_{q}^{2}}\right)\frac{\frac{1}{\sqrt{\theta}\sigma}e(z_{q}^{2})}{\overline{F}_{z}(z_{q})}.
\end{align*}
Additionally, $e(\cdot)$ is the density functions of exponential power distributions with a density generator of the form $g_{1}(u)=c_{1}\exp(-\sqrt{u})$ and $c_{1}=\frac{1}{2\sqrt{2}}$ (see Landsman and Valdez (2003)).
\end{corollary}
\noindent Proof. Letting $g_{1}(u)=c_{1}\exp(-\sqrt{2u})$, $\overline{G}_{1}(u)=c_{1}(1+\sqrt{2u})\exp(-\sqrt{2u})$ and $c_{1}=\frac{1}{2}$ (see Landsman et al. (2016)) in Theorem \ref{th.1}, we immediately obtain our result. This completes the proof of Corollary \ref{co.4}.
\section{Multivariate cases}
In this section, we consider  the multivariate TCE for mixture of elliptical distributions. To calculate it we definite shifted cumulative generator (see Landsman et al. (2016))
\begin{align}\label{(11)}
{\overline{G}}_{n-1}^{\ast}(u)=\int_{u}^{\infty}g_{n}(v+a)\mathrm{d}v,~a\geq0,~n>1,
\end{align}
with
\begin{align}\label{(12)}
{\overline{G}}_{n-1}^{\ast}(u)<\infty .
\end{align}
Here we consider $g_{n-1}^{\ast}(u)=g_{n}(u+a)$ as a density generator, if it satisfies the condition:
\begin{align}\label{(13)}
 \int_{0}^{\infty}u^{\frac{n}{2}-1}g_{n}(u+a)\mathrm{d} u<\infty,~\forall a\geq0.
 \end{align}
Let $\mathbf{Y}\sim LSME_{n}(\boldsymbol{\mu},~\mathbf{\Sigma},~\boldsymbol{\beta},~\Theta,~g_{n})$ and $\mathbf{M}=\mathbf{Y}|\Theta=\theta\sim E_{n}(\boldsymbol{\mu}+\theta\boldsymbol{\beta},~\theta\mathbf{\Sigma},~g_{n})$. Then $$\mathbf{Z}=(\theta\mathbf{\Sigma})^{-\frac{1}{2}}(\boldsymbol{M-\mu}-\theta\boldsymbol{\beta})\sim E_{n}\left(\boldsymbol{0},~\boldsymbol{I_{n}},~g_{n}\right).$$
Writing
$$\boldsymbol{\xi_{q}}=\left(\xi_{\boldsymbol{q},1},~\xi_{\boldsymbol{q},2},\cdots,\xi_{\boldsymbol{q},n}\right)^{T}=(\theta\mathbf{\Sigma})^{-\frac{1}{2}}(\boldsymbol{y_{q}-\mu}-\theta\boldsymbol{\beta}),$$
where $\boldsymbol{y_{q}}=VaR_{\boldsymbol{q}}(\boldsymbol{Y})$, and $\boldsymbol{\xi}_{\boldsymbol{q},-k}=\left(\xi_{\boldsymbol{q},1},~\xi_{\boldsymbol{q},2},\cdots,\xi_{\boldsymbol{q},k-1},~\xi_{\boldsymbol{q},k+1},\cdots,\xi_{\boldsymbol{q},n}\right)^{T}$.

To derive formula for MTCE we introduce tail function $\overline{F}_{\mathbf{Z}_{-k}}(\boldsymbol{t})$ of $(n-1)$-dimensional random vector $\mathbf{Z}_{-k}=(Z_{1},~Z_{2},\cdots,Z_{k-1},~Z_{k+1},\cdots,Z_{n})^{T}$,
$$\overline{F}_{\mathbf{Z}_{-k}}(\boldsymbol{t})=\int_{\boldsymbol{t}}^{\infty}f_{\mathbf{Z}_{-k}}(\boldsymbol{v})\mathrm{d}\boldsymbol{v},~~~\boldsymbol{v},\boldsymbol{t}\in\mathbb{R}^{n-1},~~~\mathrm{d}\boldsymbol{v}=\mathrm{d}v_{1}\mathrm{d}v_{2}\cdots\mathrm{d}v_{n},$$
with the pdf
\begin{align*}
f_{\mathbf{Z}_{-k}}(\boldsymbol{z}_{-k})&=-\frac{1}{\psi^{\ast'}(0)}\overline{G}^{\ast}_{n-1}\left\{\frac{1}{2}\boldsymbol{z}_{-k}^{T}\boldsymbol{z}_{-k}\right\}\\
&=-\frac{1}{\psi^{\ast'}(0)}\overline{G}_{n}\left\{\frac{1}{2}\boldsymbol{z}_{-k}^{T}\boldsymbol{z}_{-k}+\frac{1}{2}\xi_{\boldsymbol{q},k}^{2}\right\},~k=1,~2,\cdots,n,
\end{align*}
where
  $\psi^{\ast}(\cdot)$ is the characteristic generator corresponding to $\overline{G}^{\ast}_{n-1}$, and $\overline{G}^{\ast}_{n-1}$ as formula $(\ref{(11)})$.
\begin{theorem}\label{th.2} Let $\mathbf{Y}\sim LSME_{n}(\boldsymbol{\mu},~\mathbf{\Sigma},~\boldsymbol{\beta},~\Theta,~g_{n})$ be an $n$-dimensional location-scale mixture of elliptical random variable defined as $(\ref{(1)})$.  We suppose satisfy conditions $(\ref{(5)})$, $(\ref{(12)})$ and $(\ref{(13)})$.

 Then
\begin{align}\label{(14)}
 &MTCE_{\boldsymbol{q}}(\mathbf{Y})=\boldsymbol{\mu}+E_{\theta}\left[\theta\boldsymbol{\beta}+\sqrt{\theta}\mathbf{\Sigma}^{\frac{1}{2}}\boldsymbol{\delta_{q}}\right],
\end{align}
where
\begin{align*}
 \boldsymbol{\delta_{q}}=(\delta_{1,\boldsymbol{q}},~\delta_{2,\boldsymbol{q}},\cdots,\delta_{n,\boldsymbol{q}})^{T},
\end{align*}
with $\delta_{k,\boldsymbol{q}}=-c_{n}\psi^{\ast'}(0)\frac{\overline{F}_{\boldsymbol{z}_{-k}}(\boldsymbol{\xi}_{\boldsymbol{q},-k})}{\overline{F}_{\boldsymbol{z}}(\boldsymbol{\xi_{q}})}$ and $c_{n}=\frac{\Gamma(n/2)}{(2\pi)^{n/2}}\left[\int_{0}^{\infty}u^{\frac{n}{2}-1}g_{n}(u)\mathrm{d} u\right]^{-1}$.
\end{theorem}
\noindent Proof. Using the tower property of expectations, we obtain
\begin{align*}
MTCE_{\boldsymbol{q}}(\mathbf{Y})&=E[\mathbf{Y}|\mathbf{Y}>\boldsymbol{y_{q}}]\\
&=E_{\Theta}[E(\mathbf{Y}|\mathbf{Y}>\boldsymbol{y_{q}},\Theta)].
\end{align*}
Since
\begin{align*}
E[\mathbf{Y}|\mathbf{Y}>\boldsymbol{y_{q}},\Theta=\theta]&=E[(\mathbf{Y}|\Theta=\theta)|(\mathbf{Y}|\Theta=\theta)>\boldsymbol{y_{q}}]\\
&=E[\mathbf{M}|\mathbf{M}>\boldsymbol{y_{q}}]\\
&=MTCE_{\boldsymbol{q}}(\mathbf{M}),
\end{align*}
where $\mathbf{M}\sim E_{n}(\boldsymbol{\mu}+\theta\boldsymbol{\beta},~\theta\mathbf{\Sigma},~g_{n})$, and the second equality we have used (\ref{(2)}).
Using Theorem 1 in Landsman et al. (2016), we obtain $(\ref{(14)})$, which completes the proof of Theorem \ref{th.2}.\\
$\mathbf{Remark~4.1.}$ We remark that Theorem 1 in Landsman et al. (2016), which corresponding the result of $MTCE_{\boldsymbol{q}}(\mathbf{Y}|\Theta)$ with $\boldsymbol{q}=(q,~q,\cdots,q)^{T}$, is a special case of our Theorem 4.1.
\begin{corollary}\label{co.5}
 Suppose $\mathbf{Y}\sim LSMN_{n}(\boldsymbol{\mu},~\mathbf{\Sigma},~\boldsymbol{\beta},~\Theta)$ is an $n$-variate location-scale mixture of normal random variable defined as $(\ref{(20)})$. Under conditions in Theorem \ref{th.2}, we obtain the MTCE for location-scale mixture of normal distributions. Its' form is the same as $(\ref{(14)})$,
where
\begin{align*}
\delta_{k,\boldsymbol{q}}=-c_{n}\psi^{\ast'}(0)\frac{\overline{F}_{\boldsymbol{z}_{-k}}(\boldsymbol{\xi}_{\boldsymbol{q},-k})}{\overline{F}_{\boldsymbol{z}}(\boldsymbol{\xi_{q}})}=\phi_{1}(\xi_{k,\boldsymbol{q}})\frac{\overline{\Phi}_{\boldsymbol{z}_{-k}}(\boldsymbol{\xi}_{\boldsymbol{q},-k})}{\overline{\Phi}_{\boldsymbol{z}}(\boldsymbol{\xi_{q}})}.
\end{align*}
 Additionally, $\phi_{n}(\cdot)$ and $\Phi_{n}(\cdot)$ denote the density and distribution functions of normal distributions, respectively.
\end{corollary}
\noindent Proof. Letting the density generator $\overline{G}_{n}(u)=g_{n}(u)=\phi_{n}(u)=c_{n}e^{-u}$, $c_{n}=(2\pi)^{-\frac{n}{2}}$ and $$\psi^{\ast'}(0)=-(2\pi)^{\frac{n}{2}}\phi_{1}(\boldsymbol{\xi}_{\boldsymbol{q},k})$$  in Theorem \ref{th.2}, we directly obtain our result. This completes the proof of Corollary \ref{co.5}.
\begin{corollary}\label{co.6}
 Suppose that $\mathbf{Y}\sim LSMSt_{n}(\boldsymbol{\mu},~\mathbf{\Sigma},~\boldsymbol{\beta},~\Theta)$ is an $n$-variate location-scale mixture of student-$t$ random vector defined as $(\ref{(21)})$. Under conditions in Theorem \ref{th.2}, we obtain the MTCE for location-scale mixture of Student-$t$ distributions. Its' form is the same as $(\ref{(14)})$,
where
\begin{align*}
 &\delta_{k,\boldsymbol{q}}=-c_{n}\psi^{\ast'}(0)\frac{\overline{F}_{\boldsymbol{z}_{-k}}(\boldsymbol{\xi}_{\boldsymbol{q},-k})}{\overline{F}_{\boldsymbol{z}}(\boldsymbol{\xi_{q}})}\\
 &=\frac{\Gamma(\frac{m-1}{2})m}{2\Gamma(\frac{m}{2})\sqrt{\pi(m-1)}}\left(\frac{m-1}{m}\right)^{\frac{n}{2}}\left(1+\frac{\xi_{\boldsymbol{q},k}^{2}}{m}\right)^{-\frac{(m+n-2)}{2}}\frac{\overline{T}_{m-1,n-1}(\boldsymbol{\xi}_{\boldsymbol{q},-k};\boldsymbol{0},\Delta_{k})}{\overline{T}_{m,n}(\boldsymbol{\xi_{q}};\boldsymbol{0},\mathbf{I}_{n})},
\end{align*}
and
\begin{align*}
\Delta_{k}=\left(\frac{m(1+\frac{\xi_{\boldsymbol{q},k}^{2}}{m})}{m-1}\right)\mathbf{I}_{n-1},
\end{align*}
$\mathbf{I}_{k}$ $(k=n-1~or~n)$ is a $k$-dimensional identity matrix.
In addition, $T_{m-1,n-1}(\boldsymbol{\xi}_{\boldsymbol{q},-k};0,\Delta_{k})$ and $T_{m,n}(\boldsymbol{\xi_{q}};0,I_{n})$ are distribution functions of Student-$t$ distributions (see Landsman et al. 2016).
\end{corollary}
\noindent Proof. Letting $g_{n}(u)=c_{n}(1+\frac{2u}{m})^{-(m+n)/2}$, $\overline{G}_{n}(u)=c_{n}\frac{m}{m+n-2}(1+\frac{2u}{m})^{-(m+n-2)/2}$, $c_{n}=\frac{\Gamma((m+n)/2)}{\Gamma(m/2)(m\pi)^{\frac{n}{2}}}$ (see Landsman et al. (2016)) and
 \begin{align*}
  \psi^{\ast'}(0)=-\frac{\Gamma(\frac{m-1}{2})\pi^{(n-1)/2}(m-1)^{(n-1)/2}m}{\Gamma(\frac{m+n-2}{2})(m+n-2)}\left(1+\frac{\xi_{\boldsymbol{q},k}^{2}}{m}\right)^{-(m+n-2)/2}
  \end{align*}
 in Theorem \ref{th.2}, we immediately obtain our result. This completes the proof of Corollary \ref{co.6}.
\begin{corollary}\label{co.7}
 Assume $\mathbf{Y}\sim LSMLo_{n}(\boldsymbol{\mu},~\mathbf{\Sigma},~\boldsymbol{\beta},~\Theta)$ is an $n$-variate location-scale mixture of Logistic random vector defined as $(\ref{(22)})$. Under conditions in Theorem \ref{th.2}, we obtain the MTCE for location-scale mixture of Logistic distributions. Its' form is the same as $(\ref{(14)})$,
where
\begin{align*}
 &\delta_{k,\boldsymbol{q}}=-c_{n}\psi^{\ast'}(0)\frac{\overline{F}_{\boldsymbol{z}_{-k}}(\boldsymbol{\xi}_{\boldsymbol{q},-k})}{\overline{F}_{\boldsymbol{z}}(\boldsymbol{\xi_{q}})}\\
 &=\frac{L(-\exp(-\frac{\xi_{\boldsymbol{q},k}^{2}}{2}),~\frac{n-1}{2},~1)\exp(-\frac{\xi_{\boldsymbol{q},k}^{2}}{2})}{\sqrt{2\pi}\left[\sum_{i=0}^{\infty}(-1)^{i-1}i^{1-n/2}\right]}\frac{\overline{F}_{\boldsymbol{z}_{-k}}(\boldsymbol{\xi}_{\boldsymbol{q},-k})}{\overline{F}_{\boldsymbol{z}}(\boldsymbol{\xi_{q}})},
\end{align*}
and pdf of $~\mathbf{Z}_{-k}$:
\begin{align*}
f_{\mathbf{Z}_{-k}}(\boldsymbol{t})=-\frac{1}{\psi^{\ast'}(0)}\frac{\exp(-\frac{\boldsymbol{t}^{T}\boldsymbol{t}}{2}-\frac{\xi_{\boldsymbol{q},k}^{2}}{2})}{1+\exp(-\frac{\boldsymbol{t}^{T}\boldsymbol{t}}{2}-\frac{\xi_{\boldsymbol{q},k}^{2}}{2})},~k=1,~2,\cdots,n,~\boldsymbol{t}\in\mathbb{R}^{n-1},
\end{align*}
and \begin{align}\label{(b1)}
  \nonumber\psi^{\ast'}(0)&=-\frac{(2\pi)^{\frac{n-1}{2}}}{\Gamma(\frac{n-1}{2})}\left[\int_{0}^{\infty}\frac{t^{(n-3)/2}\exp(-t-\frac{\xi_{\boldsymbol{q},k}^{2}}{2})}{1+\exp(-t-\frac{\xi_{\boldsymbol{q},k}^{2}}{2})}\mathrm{d}t\right]\\
  &=-\frac{(2\pi)^{(n-1)/2}L(-\exp(-\frac{\xi_{\boldsymbol{q},k}^{2}}{2}),~\frac{n-1}{2},~1)}{\exp(\frac{\xi_{\boldsymbol{q},k}^{2}}{2})}.
  \end{align}
Additionally, $L(\cdot)$ is the well known Lerch zeta function (see Lin and Srivastava (2004)).
\end{corollary}
\noindent Proof. Letting $g_{n}(u)=c_{n}\frac{\exp(-u)}{[1+\exp(-u)]^{2}}$, $\overline{G}_{n}(u)=c_{n}\frac{\exp(-u)}{1+\exp(-u)}$, $$c_{n}=(2\pi)^{-n/2}\left[\sum_{i=0}^{\infty}(-1)^{i-1}i^{1-n/2}\right]^{-1}$$ and
 formula $(\ref{(b1)})$
 in Theorem \ref{th.2}, we directly obtain our result. This completes the proof of Corollary \ref{co.7}.
\begin{corollary}\label{co.8}
 Assume that $\mathbf{Y}\sim LSMLa_{n}(\boldsymbol{\mu},~\mathbf{\Sigma},~\boldsymbol{\beta},~\Theta)$ be an $n$-variate location-scale mixture of Laplace random vector defined as $(\ref{(22)})$. Under conditions in Theorem \ref{th.2}, we obtain the MTCE for location-scale mixture of Laplace distributions. Its' form is the same as $(\ref{(14)})$,
where
\begin{align*}
 &\delta_{k,\boldsymbol{q}}=-c_{n}\psi^{\ast'}(0)\frac{\overline{F}_{\boldsymbol{z}_{-k}}(\boldsymbol{\xi}_{\boldsymbol{q},-k})}{\overline{F}_{\boldsymbol{z}}(\boldsymbol{\xi_{q}})}\\
 &=-\frac{\Gamma(n/2)\psi^{\ast'}(0)}{2\pi^{n/2}\Gamma(n)}\frac{\overline{F}_{\boldsymbol{z}_{-k}}(\boldsymbol{\xi}_{\boldsymbol{q},-k})}{\overline{F}_{\boldsymbol{z}}(\boldsymbol{\xi_{q}})},
\end{align*}
and pdf of $\mathbf{Z}_{-k}$:
\begin{align*}
f_{\mathbf{Z}_{-k}}(\boldsymbol{t})=-\frac{1}{\psi^{\ast'}(0)}\left(1+\sqrt{\boldsymbol{t}^{T}\boldsymbol{t}+\xi_{\boldsymbol{q},k}^{2}}\right)\exp\left\{-\sqrt{\boldsymbol{t}^{T}\boldsymbol{t}+\xi_{\boldsymbol{q},k}^{2}}\right\},~k=1,~2,\cdots,n,~\boldsymbol{t}\in\mathbb{R}^{n-1},
\end{align*}
and
\begin{align}\label{(b2)}
\psi^{\ast'}(0)=-\frac{(2\pi)^{(n-1)/2}}{\Gamma\left(\frac{n-1}{2}\right)}\left[\int_{0}^{\infty}t^{\frac{n-3}{2}}\left(1+\sqrt{t+\xi_{\boldsymbol{q},k}^{2}}\right)\exp\left\{-\sqrt{t+\xi_{\boldsymbol{q},k}^{2}}\right\}\mathrm{d}t\right].
\end{align}
\end{corollary}
\noindent Proof. Letting $g_{n}(u)=c_{n}\exp(-\sqrt{2u})$, $\overline{G}_{n}(u)=c_{n}(1+\sqrt{2u})\exp(-\sqrt{2u})$, $c_{n}=\frac{\Gamma(n/2)}{2\pi^{n/2}\Gamma(n)}$ (see Landsman et al. 2016) and formula $(\ref{(b2)})$
in Theorem \ref{th.2}, we immediately obtain our result. This completes the proof of Corollary \ref{co.8}.
\section{Portfolio risk decomposition with TCE}
Let $\mathbf{Y}=(Y_{1},~Y_{2},\cdots,Y_{n})^{T}\sim LSME_{n}(\boldsymbol{\mu},~\mathbf{\Sigma},~\boldsymbol{\beta},~\Theta,~g_{n})$, $\boldsymbol{e}=(1,~1,\cdots,1)^{T}$ is an $n\times1$ vector whose elements are all equal to $1$. We define
\begin{align}\label{(15)}
  S=\sum_{j=1}^{n}Y_{j}=\boldsymbol{e}^{T}\mathbf{Y}=\boldsymbol{e}^{T}\boldsymbol{\mu}+\Theta\boldsymbol{e}^{T}\boldsymbol{\beta}+\Theta^{\frac{1}{2}}\boldsymbol{e}^{T}\mathbf{\Sigma}^{\frac{1}{2}}\mathbf{X},
\end{align}
which is the sum of mixture of elliptical risks.\\
$\mathbf{Proposition~5.1.}$
Under the conditions $(\ref{(9)})$ and $(\ref{(15)})$, the TCE of $S$ can be expressed as
\begin{align}\label{(17)}
 &TCE_{S}(s_{q})=\mu_{S}+E_{\theta}\left[\theta\beta_{S}+\delta_{S}(\sqrt{\theta}\sigma_{S})^{2}\right],
\end{align}
where
\begin{align*}
 \delta_{S}=\frac{\frac{1}{\sqrt{\theta}\sigma_{S}}\overline{G}_{1}(\frac{1}{2}z_{q}^{2})}{\overline{F}_{Z}(z_{q})},
\end{align*}
with $Z\sim E_{1}(0,~1,~g_{1})$ and $z_{q}=\frac{s_{q}-\mu_{S}-\theta\beta_{S}}{\sqrt{\theta}\sigma_{S}}$.\\
\noindent Proof. Let $L=\boldsymbol{e}^{T}\mathbf{\Sigma}^{\frac{1}{2}}\mathbf{X}$. Due to $(\ref{(6)})$, we get
$L\sim E_{1}(0,~\sigma_{S}^{2},~g_{1})$ with $\sigma_{S}^{2}=\boldsymbol{e}^{T}\mathbf{\Sigma}\boldsymbol{e}$. So that \\$L'\sim E_{1}(0,~1,~g_{1})$ with $L'=\frac{L}{\sigma_{S}}$. Therefore,
\begin{align}\label{(16)}
S=\mu_{S}+\Theta\beta_{S}+\Theta^{\frac{1}{2}}\sigma_{L}L'\sim LSME_{1}(\mu_{S},~\sigma_{S}^{2},~\beta_{S},~\Theta,~g_{1}),
\end{align}
with $\mu_{S}=\boldsymbol{e}^{T}\boldsymbol{\mu}$ and $\beta_{S}=\boldsymbol{e}^{T}\boldsymbol{\beta}$.\\
By using Theorem \ref{th.1}, we obtain $(\ref{(17)})$, which completes the proof of Proposition 5.1.
\begin{lemma}\label{le.1}
Let $\mathbf{Y}=(Y_{1},~Y_{2},\cdots,Y_{n})^{T}\sim LSME_{n}(\boldsymbol{\mu},~\mathbf{\Sigma},~\boldsymbol{\beta},~\Theta,~g_{n})$ as $(\ref{(1)})$. Then the vector \\$Y_{k,S}=(Y_{k},~S)^{T}$, $(1\leq k\leq n)$ has a mixture of elliptical distribution, namely, $$Y_{k,S}\sim LSME_{2}(\boldsymbol{\mu}_{k,S},~\mathbf{\Sigma}_{k,S},~\boldsymbol{\beta}_{k,S},~\Theta,~g_{2}),$$
 where $\boldsymbol{\mu}_{k,S}=(\mu_{k},~\boldsymbol{e}^{T}\boldsymbol{\mu})^{T}=(\mu_{k},~\sum_{i=1}^{n}\mu_{i})^{T}$,
\begin{align*}
\mathbf{\Sigma}_{k,S}=\left(\begin{array}{ccccccccccc}
\sigma_{k}^{2}&\sigma_{k,S}\\
\Sigma_{k,S}&\sigma_{S}^{2}
\end{array}
\right)
\end{align*}
and $\boldsymbol{\beta}_{k,S}=(\beta_{k},~\beta_{S})^{T}$, and $\sigma_{k}^{2}=\sigma_{k,k}$, $\sigma_{k,S}=\sum_{i=1}^{n}\sigma_{k,i}$, $\beta_{S}=\boldsymbol{e}^{T}\boldsymbol{\beta}=\sum_{i=1}^{n}\beta_{i}$,\\
$\sigma_{S}^{2}=\boldsymbol{e}^{T}\mathbf{\Sigma} \boldsymbol{e}=\sum_{i,j=1}^{n}\sigma_{i,j}.$
\end{lemma}
\noindent Proof. Write $\mathbf{P}=(P_{1},~P_{2},\cdots,P_{n})^{T}=\mathbf{\Sigma}^{\frac{1}{2}}\mathbf{X}$. Due to $(\ref{(6)})$, we know
$\mathbf{P}\sim E_{n}(\boldsymbol{0},~\mathbf{\Sigma},~g_{n}).$\\
From $(\ref{(1)})$, we get $Y_{k}=\mu_{k}+\Theta\beta_{k}+\Theta^{\frac{1}{2}}P_{k}$, $1\leq k\leq n$.
According to $(\ref{(16)})$, we know $S=\mu_{S}+\Theta\beta_{S}+\Theta^{\frac{1}{2}}L$ with $L\sim E_{1}(0,~\sigma_{S}^{2},~g_{1})$. So that $Y_{k,S}=\boldsymbol{\mu}_{k,S}+\Theta\boldsymbol{\beta}_{k,S}+\Theta^{\frac{1}{2}}(P_{k},~L)^{T}.$\\
By Lemma 1 in Landsman and Valdez (2003), we obtain $(P_{k},~L)^{T}\sim E_{2}(\boldsymbol{\mu}_{k,S},~\mathbf{\Sigma}_{k,S},~g_{2})$. Therefore, $Y_{k,S}\sim LSME_{2}(\boldsymbol{\mu}_{k,S},~\mathbf{\Sigma}_{k,S},~\boldsymbol{\beta}_{k,S},~\Theta,~g_{2})$. This completes the proof of Lemma \ref{le.1}.
\begin{lemma}\label{le.2}
Let $\mathbf{Y}=(Y_{1},~Y_{2})^{T}\sim LSME_{2}(\boldsymbol{\mu},~\mathbf{\Sigma},~\boldsymbol{\beta},~\Theta,~g_{2})$. Assume that
condition $(\ref{(9)})$ holds. Then
\begin{align}\label{(18)}
TCE_{Y_{1}|Y_{2}}(y_{q})=\mu_{1}+E_{\theta}[\theta\beta_{1}+\delta_{2}\theta\sigma_{1}\sigma_{2}\rho_{1,2}],
\end{align}
where
$$\delta_{2}=\frac{\frac{1}{\sqrt{\theta}\sigma_{2}}\overline{G}(\frac{1}{2}z_{2,q}^{2})}{\overline{F}_{z}(z_{2,q})},$$
$\rho_{1,2}=\frac{\sigma_{1,2}}{\sigma_{1}\sigma_{2}}$, $\sigma_{1}=\sqrt{\sigma_{1,1}}$, $\sigma_{2}=\sqrt{\sigma_{2,2}}$ and $z_{2,q}=\frac{y_{q}-\mu_{2}-\theta\beta_{2}}{\sqrt{\theta}\sigma_{2}}$.
\end{lemma}
\noindent Proof.  Using the tower property of expectations, we have
\begin{align*}
TCE_{Y_{1}|Y_{2}}(y_{q})&=E[Y_{1}|Y_{2}>y_{q}]\\
&=E_{\theta}[E[Y_{1}|Y_{2}>y_{q},\Theta=\theta]]\\
&=E_{\theta}[E[Q_{1}|Q_{2}>y_{q}]],
\end{align*}
where $(Q_{1},~Q_{2})^{T}=\mathbf{Y}|\Theta=\theta\sim E_{2}(\boldsymbol{\mu}+\theta\boldsymbol{\beta},~\theta\mathbf{\Sigma},~g_{2})$, and the third equality we have used (\ref{(2)}).\\
By Lemma 2 in Landsman and Valdez (2003), we can obtain $(\ref{(18)})$. This completes the proof of Lemma \ref{le.2}.

We use the above two Lemmas to give portfolio risk decomposition with TCE as follows.
\begin{theorem}\label{th.3} Let $\mathbf{Y}=(Y_{1},~Y_{2},\cdots,Y_{n})^{T}\sim LSME_{n}(\boldsymbol{\mu},~\mathbf{\Sigma},~\boldsymbol{\beta},~\Theta,~g_{n})$ be an $n$-dimensinal location-scale mixture of elliptical random vector defined as $(\ref{(1)})$.  We suppose condition $(\ref{(9)})$ holds, and let $S=\sum_{i=1}^{n}Y_{i}$. Then the contribution of risk $Y_{k},~1\leq k\leq n$, to the total TCE can be given by
\begin{align}\label{(19)}
TCE_{Y_{k}|S}(s_{q})=\mu_{k}+E_{\theta}[\theta\beta_{k}+\delta_{S}\theta\sigma_{k}\sigma_{S}\rho_{k,S}],
\end{align}
where $\rho_{k,S}=\frac{\sigma_{k,S}}{\sigma_{k}\sigma_{S}}$ and $\delta_{S}$ is the same as in Proposition 5.1.
\end{theorem}
\noindent Proof. By Lemma \ref{le.1}, we know
$Y_{k,S}=(Y_{k},~S)^{T}\sim LSME_{2}(\boldsymbol{\mu}_{k,S},~\mathbf{\Sigma}_{k,S},~\boldsymbol{\beta}_{k,S},~\Theta,~g_{2})$, $(1\leq k\leq n)$.
Let $\mathbf{Y}$ subject to $\mathbf{Y}_{k,S}$ in Lemma \ref{le.2}, we can immediately obtain $(\ref{(19)})$. This completes the proof of Theorem \ref{th.3}.\\
$\mathbf{Remark~5.1.}$  Letting the density generator $\overline{G}_{1}(u)=g_{1}(u)=\phi_{1}(u)=(2\pi)^{-\frac{1}{2}}e^{-u}$ in Theorem \ref{th.3}, we obtain the portfolio risk decomposition with TCE for location-scale mixture of normal distributions. Its' form is the same as $(\ref{(19)})$,
where
\begin{align*}
 \delta_{S}=\frac{\frac{1}{\sqrt{\theta}\sigma_{S}}\phi_{1}(\frac{1}{2}z_{q}^{2})}{1-\Phi_{1}(z_{q})}.
\end{align*}
 Additionally, $\phi_{1}(\cdot)$ and $\Phi_{1}(\cdot)$ denote the density and distribution functions of normal distributions.\\
$\mathbf{Remark~5.2.}$ Letting $g_{1}(u)=c_{1}(1+\frac{2u}{m})^{-(m+1)/2}$, $\overline{G}_{1}(u)=c_{1}\frac{m}{m-1}(1+\frac{2u}{m})^{-(m-1)/2}$ and $c_{1}=\frac{\Gamma((m+1)/2)}{\Gamma(m/2)(m\pi)^{\frac{1}{2}}}$ (see Landsman et al. (2016)) in Theorem \ref{th.3}, we obtain the portfolio risk decomposition with TCE for location-scale mixture of Student-$t$ distributions. Its' form is the same as $(\ref{(19)})$,
where
\begin{align*}
 \delta_{S}=\frac{\frac{1}{\sqrt{\theta}\sigma_{S}}\overline{G}_{1}(\frac{1}{2}z_{q}^{2})}{\overline{F}_{Z}(z_{q})}=\frac{\frac{1}{\sqrt{\theta}\sigma_{S}}c_{1}\frac{m}{m-1}(1+\frac{z_{q}^{2}}{m})^{-(m-1)/2}}{\overline{F}_{Z}(z_{q})}=\frac{\frac{1}{\sqrt{\theta}\sigma_{S}}t_{m,1}(z_{q};0,1)}{\overline{T}_{m,1}(z_{q};0,1)}.
\end{align*}
In addition, $t_{m,1}(z_{q};0,1)$ and $T_{m,1}(z_{q};0,1)$ are the density and distribution functions of Student-$t$ distributions, respectively (see Landsman et al. (2016)).\\
$\mathbf{Remark~5.3.}$ Letting $g_{1}(u)=c_{1}\frac{\exp(-u)}{[1+\exp(-u)]^{2}}$, $\overline{G}_{1}(u)=c_{1}\frac{\exp(-u)}{1+\exp(-u)}$ and $c_{1}=\frac{1}{2}$ (see Landsman and Valdez (2003)) in Theorem \ref{th.3}, we obtain the portfolio risk decomposition with TCE for location-scale mixture of Logistic distributions. Its' form is the same as $(\ref{(19)})$,
where
\begin{align*}
 \delta_{S}=\frac{\frac{1}{\sqrt{\theta}\sigma_{S}}\overline{G}_{1}(\frac{1}{2}z_{q}^{2})}{\overline{F}_{Z}(z_{q})}=\frac{\frac{1}{\sqrt{\theta}\sigma_{S}}c_{1}\frac{\exp(-\frac{1}{2}z_{q}^{2})}{1+\exp(-\frac{1}{2}z_{q}^{2})}}{\overline{F}_{Z}(z_{q})}=\left[\frac{1}{2[(\sqrt{2\pi})^{-1}+\phi(z_{q})]}\right]\frac{\frac{1}{\sqrt{\theta}\sigma_{S}}\phi(z_{q})}{\overline{F}_{z}(z_{q})}.
\end{align*}
In addition, $\phi(\cdot)$ is the density functions of normal distributions (see Landsman and Valdez (2003)).\\
$\mathbf{Remark~5.4.}$ Letting $g_{1}(u)=c_{1}\exp(-\sqrt{2u})$, $\overline{G}_{1}(u)=c_{1}(1+\sqrt{2u})\exp(-\sqrt{2u})$ and $c_{1}=\frac{1}{2}$ (see Landsman et al. (2016)) in Theorem \ref{th.3}, we obtain the portfolio risk decomposition with TCE for location-scale mixture of Laplace distributions. Its' form is the same as $(\ref{(19)})$,
where
\begin{align*}
 \delta_{S}=\frac{\frac{1}{\sqrt{\theta}\sigma_{S}}\overline{G}_{1}(\frac{1}{2}z_{q}^{2})}{\overline{F}_{Z}(z_{q})}=\frac{\frac{1}{\sqrt{\theta}\sigma_{S}}c_{1}(1+\sqrt{z_{q}^{2}})\exp(-\sqrt{z_{q}^{2}})}{\overline{F}_{Z}(z_{q})}=\sqrt{2}\left(1+\sqrt{z_{q}^{2}}\right)\frac{\frac{1}{\sqrt{\theta}\sigma_{S}}e(z_{q}^{2})}{\overline{F}_{z}(z_{q})}.
\end{align*}
Additionally, $e(\cdot)$ is the density functions of exponential power distributions with a density generator of the form $g_{1}(u)=c_{1}\exp(-\sqrt{u})$ and $c_{1}=\frac{1}{2\sqrt{2}}$ (see Landsman and Valdez (2003)).
\section{Concluding remarks}
 In this paper we consider the univariate and multivariate location-scale mixture of elliptical distribution, which is $(\mathbf{A}=\mathbf{\Sigma}^{\frac{1}{2}})$ generalization of normal mean-variance mixture distribution in Kim and Kim (2019). It has received much attention in finance and insurance applications, since this distribution not only include location-scale mixture of normal (LSMN) distributions, location-scale mixture of Student-$t$ (LSMSt) distributions, location-scale mixture of Logistic (LSMLo) distributions and location-scale mixture of Laplace (LSMLa) distributions, but also include the generalized hyperbolic distribution (GHD) and the slash distribution. The GHD is a special case of this mixture random variable with $\mathbf{X}\sim N_{n}(\boldsymbol{0},~\mathbf{I}_{n})$ and the distribution of $\Theta$ given by a generalized inverse gaussian $N^{-1}(\lambda,~\chi,~\psi)$ (see Kim and Kim (2019) for details). The GHD is an important distribution, and has a lot of applications (see Kim (2010) and Ignatieva and Landsman (2015)). Slash distribution also is a special case of this mixture random variable with $\mathbf{X}\sim N_{n}(\boldsymbol{0},~\mathbf{I}_{n})$ and $\Theta\sim BP(\eta=1,~\alpha=1,~\beta=q/2)$. Here $BP(\cdot)$ is the $3$-parameter beta prime (BP) or inverted beta distribution (see Kim and Kim (2019) for details). This distribution has been discussed in many literatures (see Gneiting (1997), Gen\c{c} (2007) and Wang and Genton (2006)). We also consider univariate TCE, multivariate TCE and portfolio risk decomposition with TCE for location-scale mixture of elliptical distribution. As special cases, we provided univariate TCE, multivariate TCE and portfolio risk decomposition with TCE for LSMN, LSMSt, LSMLo and LSMLa distributions.
\section*{Acknowledgments}
\noindent The research was supported by the National Natural Science Foundation of China (No.11571198, 11701319)




\end{document}